\documentclass[11pt,leqno]{article} 
\usepackage{graphics}
\newtheorem{thm}{Theorem}[section]
\newtheorem{lma}{Lemma}[section]
\newtheorem{cor}{Corollary}

\newcommand{\beqa}{\begin{eqnarray}}
\newcommand{\eeqa}{\end{eqnarray}}

\newcommand{\pf}{\noindent {\bf Proof:} $\s$ }
\newcommand{\epf}{ \hfill$\diamondsuit$ \medskip}

\newcommand{\ds}{\displaystyle}
\newcommand{\beq}{\begin{equation}}
\newcommand{\eeq}{\end{equation}}
\newcommand{\lbl}{\label}
\newcommand{\s}{\; \;}

\newcommand{\la}{\lambda}

\newcommand{\ra}{\rightarrow}

\newcommand{\p}{\varphi}
\title{Curves of equiharmonic solutions, and problems at resonance}

\author{
Philip Korman   \\ 
Department of Mathematical Sciences \\ 
University of Cincinnati \\ 
Cincinnati Ohio 45221-0025 \\
}

\date{}

\begin{document}

\maketitle
\begin{abstract} 
We consider the semilinear Dirichlet problem
\[
\Delta u+kg(u)=\mu _1 \p _1+\cdots +\mu _n \p _n+e(x) \s \mbox{for $x \in \Omega$}, \s u=0 \s \mbox{on $\partial \Omega$},
\]
where 
$\p _k$ is the $k$-th eigenfunction of the Laplacian on $\Omega$ and $e(x) \perp \p _k$, $k=1, \ldots,n$. Write the solution in the form $u(x)= \Sigma _{i=1}^n  \xi _i \p _i+U(x)$, with $ U \perp \p _k$, $k=1, \ldots,n$. Starting with $k=0$, when the problem is linear, we continue the solution in $k$ by keeping $\xi =(\xi _1, \ldots,\xi _n)$ fixed, but allowing for $\mu =(\mu _1, \ldots,\mu _n)$ to vary. Studying  the map $\xi \ra \mu$ provides us with the existence and multiplicity results for the above problem. We apply our results to problems at resonance, at both the  principal and higher eigenvalues. Our approach is suitable for numerical calculations, which we  implement, illustrating our results.
 \end{abstract}

\begin{flushleft}
\end{flushleft}

\begin{flushleft}
\end{flushleft}

\section{Introduction}
\setcounter{equation}{0}
\setcounter{thm}{0}
\setcounter{lma}{0}

We study existence and multiplicity of solutions for a semilinear problem 
\beqa
\lbl{0}
& \Delta u+kg(u)=f(x) \s \mbox{for $x \in \Omega$} \,, \\ \nonumber
& u=0 \s \mbox{on $\partial \Omega$}
\eeqa
on a smooth bounded domain $\Omega \subset R^m$. Here the functions $f(x) \in L^2(\Omega)$ and $g(u) \in C^1(R)$ are given, $k$ is a parameter. We approach this problem by continuation in $k$.  When $k=0$ the problem is linear. It has a unique solution, as can be seen by using Fourier series of the form $u(x)=\Sigma _{j=1}^{\infty} u_j \p _j$, where 
$\p _j$ is the $j$-th eigenfunction of the Dirichlet Laplacian on $\Omega$, with $\int_\Omega \p _j ^2 \, dx=1$, and $\la _j$ is the corresponding eigenvalue. We now continue the solution in $k$, looking for a solution pair $(k,u)$, or  $u=u(x,k)$. At a generic point $(k,u)$ the implicit function theorem applies, allowing the continuation in $k$. These are the {\em regular} points, where the corresponding linearized problem has only the trivial solution. So until a {\em singular} point is encountered, we have a  solution curve $u=u(x,k)$. At a singular point  practically anything imaginable might happen. At some singular points the M.G. Crandall and P.H. Rabinowitz bifurcation theorem \cite{CR} applies, giving us a curve of solutions through a singular point. But even in this favorable situation there is a possibility that solution curve will ``turn back" in $k$.
\medskip 

In \cite{K1} we have presented a  way to continue solutions forward in $k$, which can take us through any singular point. We describe it next.  If a solution $u(x)$ is given by its Fourier series  $u(x)=\Sigma _{j=1}^{\infty} \xi _j \p _j$, we call $U_n=(\xi _1,\xi _2,\ldots, \xi _n)$ the {\em $n$-signature} of solution, or just {\em signature} for short. We also represent  $f(x)$ by its Fourier series, and rewrite the problem (\ref{0}) as 
\beqa
\lbl{0.1}
& \Delta u+kg(u)=\mu^0 _1 \p _1+\cdots +\mu^0 _n \p _n+e(x) \s \mbox{for $x \in \Omega$}, \\ \nonumber
& u=0 \s \mbox{on $\partial \Omega$}
\eeqa
with $\mu^0 _j=\int _\Omega f \p _j \, dx$, and $e(x)$ is the projection of $f(x)$ onto  the orthogonal complement to $\p _1, \ldots, \p_n$. Let us now constrain ourselves to hold the signature  $U_n$ fixed (when continuing in $k$), and in return allow for $\mu _1, \ldots, \mu_n$ to vary. I.e., we are looking for $(u, \mu _1, \ldots, \mu_n)$ as a function of $k$, with $U_n$ fixed, solving
\beq
\lbl{0.1a}
\s\s\s \Delta u+kg(u)=\mu _1 \p _1+\cdots +\mu _n \p _n+e(x) \s \mbox{for $x \in \Omega$}, \s u=0 \s \mbox{on $\partial \Omega$} \,,
\eeq
\[
\int _\Omega u \p _i \, dx=\xi _i, \s i=1, \dots, n \,.
\]
It turned out that we can continue forward in $k$ this way, so long as 
\beq
\lbl{0.2}
k \max _{u \in R} g'(u)<\la _{n+1} \,.
\eeq
In the present paper we present a much simplified proof of this result, and generalize it for the case of $(i,n)$ signatures (defined below). Then, we present two new applications.
\medskip

So suppose the condition (\ref{0.2}) holds, and we wish to solve the problem (\ref{0.1}) at some $k=k_0$. We travel in $k$, from $k=0$ to $k=k_0$, on a curve of fixed signature $U_n=(\xi _1,\xi _2,\ldots, \xi _n)$, obtaining a solution $(u, \mu _1, \ldots, \mu_n)$ of (\ref{0.1a}). The right hand side of (\ref{0.1a}) has the first $n$ harmonics different (in general) from the ones we want in (\ref{0.1}). We now vary $U_n$. The question is: can we choose $U_n$ to obtain the desired $\mu _1=\mu^0 _1, \ldots, \mu_n=\mu^0 _n$, and if so, in how many ways? This corresponds to the existence and multiplicity questions for the original problem (\ref{0}).
In \cite{K1} we obtained  this way a unified approach  
to the well known  results of E.M. Landesman and A.C. Lazer \cite{L}, A. Ambrosetti and G.  Prodi \cite{AP}, M. S. Berger and E.  Podolak \cite{BP},  H. Amann and P. Hess \cite{AH} and D.G. de Figueiredo and W.-M.  Ni \cite{FN}. We also provided some new results on ``jumping nonlinearities", and on symmetry breaking.
\medskip

Our main new application in the present paper is to unbounded perturbations at resonance, which we describe next. For the problem
\[
\Delta u +\la _1 u+g(u)=e(x) \s \mbox{on $\Omega$}, \s u=0 \s \mbox{on $\partial \Omega$} \,,
\]
with a {\em bounded} $g(u)$, satisfying $ug(u) \geq 0$ for all $u \in R$, and  $e(x) \in L^{\infty} (\Omega)$ satisfying $\int _\Omega e(x) \p _1(x) \, dx=0$, D.G. de Figueiredo and W.-M.  Ni \cite{FN} have proved the existence of solutions. R. Iannacci, M.N. Nkashama and J.R. Ward \cite{I} generalized this result
to unbounded $g(u)$ satisfying $g'(u) \leq \gamma <\la _2-\la _1$ (they can also treat the case $\gamma =\la _2-\la _1$ under an additional condition).
We consider a more general problem 
\[
\Delta u +\la _1 u+g(u)=\mu _1 \p_1 +e(x) \s \mbox{on $\Omega$}, \s u=0 \s \mbox{on $\partial \Omega$} \,,
\]
with $g(u)$ and $e(x)$ satisfying the same conditions. Writing $u=\xi _1 \p _1+U$, we show that there exists a continuous curve of solutions $(u,\mu _1)(\xi _1)$, and all solutions lie on this curve. Moreover $\mu _1(\xi _1)> 0$ ($<0$) for $\xi _1>0$ ($<0$) and large. By continuity, $\mu _1(\xi ^0_1)=0$ at some $ \xi ^0_1$. We see that the existence result of R. Iannacci et al \cite{I} corresponds to just one point on this solution  curve.
\medskip

Our second application is to resonance at higher eigenvalues, where we operate with multiple harmonics. We obtain an extension of D.G. de Figueiredo and W.-M.  Ni's \cite{FN} result to any simple $\la _k$. 
\medskip

Our approach in the present paper is well suited for numerical computations. We describe the implementation of the numerical computations, and use them to give numerical examples for our results.

\section{Preliminary results}
\setcounter{equation}{0}
\setcounter{thm}{0}
\setcounter{lma}{0}

Recall that on a smooth bounded domain $\Omega \subset R^m$ the eigenvalue problem
\[
\Delta u +\la u=0 \s \mbox{on $\Omega$}, \s u=0 \s \mbox{on $\partial \Omega$}
\]
has an infinite sequence of eigenvalues $0<\la _1<\la _2 \leq \la _3\leq \ldots \ra \infty$, where we repeat each eigenvalue according to its multiplicity, and the corresponding eigenfunctions we denote $\varphi _k$. These eigenfunctions $\varphi _k$ form an orthogonal basis of $L^2(\Omega)$, i.e., any $f(x) \in L^2(\Omega)$ can be written as $f(x)=\Sigma _{k=1}^{\infty} a_k \p _k$, with the series convergent in $L^2(\Omega)$, see e.g., L. Evans \cite{E}. We  normalize $||\varphi _k||_{L^2(\Omega)}=1$, for all $k$.

\begin{lma}\lbl{lma:1}
Assume that $u(x) \in L^2(\Omega)$, and $u(x)=\Sigma _{k=n+1}^{\infty} \xi _k \p _k$. Then
\[
\int_\Omega |\nabla u|^2 \, dx  \geq \la _{n+1} \int _\Omega u^2 \, dx.
\]
\end{lma}

\pf
Since $u(x)$ is orthogonal to $\p _1, \, \ldots, \p _n$, the proof follows by the variational characterization of $\la _{n+1}$.
\epf

In the following linear problem the function $a(x)$ is given, while $\mu _1, \ldots, \mu _n$, and $w(x)$ are unknown.

\begin{lma}\lbl{lma:3}
Consider the problem
\beqa
\lbl{9}
& \Delta w+a(x)w=\mu_1 \p _1+ \cdots +\mu _n \p _n, \s  \, \mbox{for $x \in \Omega$}, \\ \nonumber
& w=0 \s \mbox{on $\partial \Omega$}, \\ \nonumber
& \int _\Omega w \p _1 \, dx= \ldots = \int _\Omega w \p _n \, dx=0.
\eeqa
Assume that 
\beq
\lbl{10}
a(x) < \la _{n+1}, \,  \s \mbox{for all $x \in \Omega$}.
\eeq
Then the only solution of (\ref{9}) is $\mu _1 = \ldots =\mu _n=0$, and $w(x) \equiv 0$.
\end{lma}

\pf
Multiply the equation in (\ref{9}) by $w(x)$, a solution of the problem (\ref{9}), and integrate. Using Lemma \ref{lma:1} and the assumption (\ref{10}), we have
\[
\la _{n+1} \int _\Omega w^2 \, dx \leq \int_\Omega |\nabla w|^2 \, dx  =\int_\Omega a(x) w^2  \, dx < \la _{n+1} \int _\Omega w^2 \, dx.
\]
It follows that $w(x) \equiv 0$, and then
\[
0=\mu_1 \p _1+ \cdots +\mu _n \p _n \s \mbox{for $x \in \Omega$},
\]
which implies that $\mu _1 = \ldots =\mu _n=0$.
\epf

\begin{cor}\lbl{cor:1}
If one considers the problem (\ref{9}) with $\mu _1 = \ldots = \mu _n =0$, then $w(x) \equiv 0$ is the only solution of that problem.
\end{cor}

\begin{cor}\lbl{cor:2}
With $f(x) \in L^2(\Omega)$, consider the problem
\beqa \nonumber
& \Delta w+a(x)w=f(x) \s\s \mbox{for $x \in \Omega$}\,, \\ \nonumber
& w=0 \s\s \mbox{on $\partial \Omega$}, \\ \nonumber
& \int _\Omega w \p _1 \, dx= \ldots = \int _\Omega w \p _n \, dx=0.
\eeqa
Then there is a constant $c$, so that the following a priori estimate holds
\[
||w||_{H^2(\Omega)} \leq c ||f||_{L^2(\Omega)} \,.
\]
\end{cor}

\pf
An elliptic estimate gives
\[
||w||_{H^2(\Omega)} \leq c \left(||w||_{L^2(\Omega)}+||f||_{L^2(\Omega)} \right) \,.
\]
Since the corresponding homogeneous problem has only the trivial solution, the extra term on the right is removed in a standard way.
\epf

We shall also need a variation of the above lemma.

\begin{lma}\lbl{lma:4}
Consider the problem ($2 \leq i <n$)
\beqa
\lbl{11}
& \Delta w+a(x)w=\mu_i \p _i+\mu_{i+1} \p _{i+1}+ \cdots +\mu _n \p _n \s \mbox{for $x \in \Omega$}, \\ \nonumber
& w=0 \s \mbox{on $\partial \Omega$}, \\ \nonumber
& \int _\Omega w \p _i \, dx= \int _\Omega w \p _{i+1} \, dx=\ldots = \int _\Omega w \p _n \, dx=0.
\eeqa
Assume that 
\beq
\lbl{12}
\la _{i-1} <a(x) < \la _{n+1}, \,  \s \mbox{for all $x \in \Omega$} \,.
\eeq
Then the only solution of (\ref{11}) is $\mu _i = \ldots =\mu _n=0$, and $w(x) \equiv 0$.
\end{lma}

\pf
Since the harmonics from $i$-th to $n$-th are missing in the solution, we may represent $\ds w=w_1+w_2$, with $\ds w_1 \in Span \{ \p_1, \ldots, \p _{i-1} \}$, and $\ds w_2 \in Span \{ \p_{n+1}, \p _{n+2}, \ldots \}$. Multiply the equation (\ref{11}) by $w_1$, and integrate
\[
-\int _\Omega |\nabla w _1|^2 \, dx+\int _\Omega a(x) w _1^2 \, dx+\int _\Omega a(x) w _1 w_2 \, dx=0 \,.
\]
Similarly
\[
-\int _\Omega |\nabla w _2|^2 \, dx+\int _\Omega a(x) w _2^2 \, dx+\int _\Omega a(x) w _1 w_2 \, dx=0 \,.
\]
Subtracting
\beq
\lbl{14}
\int _\Omega |\nabla w _2|^2 \, dx-\int _\Omega |\nabla w _1|^2 \, dx=\int _\Omega a(x) w _2^2 \, dx-\int _\Omega a(x) w _1^2 \, dx \,.
\eeq
By the variational characterization of eigenvalues, the quantity on the left in (\ref{14}) is greater or equal to
\[
\la _{n+1} \int _\Omega  w _2^2 \, dx-\la _{i-1} \int _\Omega  w _1^2 \, dx \,,
\]
while the one of the on the right is strictly less than the above number, by our condition (\ref{12}). We have a contradiction, unless $w_1=w_2 \equiv 0$. Then $\mu _i = \ldots =\mu _n=0$.
\epf

\begin{cor}\lbl{cor:3}
If one considers the problem (\ref{11}) with $\mu _i = \ldots = \mu _n =0$, then $w(x) \equiv 0$ is the only solution of that problem. Consequently, for the problem 
\beqa \nonumber
& \Delta w+a(x)w=f(x) \s\s \mbox{for $x \in \Omega$}\,, \\ \nonumber
& w=0 \s\s \mbox{on $\partial \Omega$}, \\ \nonumber
& \int _\Omega w \p _i \, dx= \ldots = \int _\Omega w \p _n \, dx=0.
\eeqa
there is a constant $c$, so that the following a priori estimate holds
\[
||w||_{H^2(\Omega)} \leq c ||f||_{L^2(\Omega)} \,.
\]
\end{cor}

\section{Continuation of solutions}
\setcounter{equation}{0}
\setcounter{thm}{0}
\setcounter{lma}{0}

Any $f(x) \in L^2(\Omega)$ can be decomposed as $f(x)=\mu_1 \p _1+ \ldots +\mu _n \p _n+e(x)$, with $e(x)=\Sigma _{j=n+1}^{\infty} e _j \p _j$ orthogonal to $\p _1,  \ldots, \p _n$. We consider a boundary value problem
\beqa
\lbl{2}
& \Delta u+kg(u)=\mu_1 \p _1+ \ldots +\mu _n \p _n+e(x) \s \mbox{for $x \in \Omega$}, \\ \nonumber
& u=0 \s \mbox{on $\partial \Omega$}.
\eeqa
Here $k \geq 0$ is a constant, and $g(u) \in C^1(R) $ is assumed to satisfy
\beq
\lbl{4}
g(u)= 
\gamma u+b(u) \,,
\eeq
with a real constant $\gamma $,  and $b(u)$ bounded for all $u \in R$, and also
\beq
\lbl{3}
g'(u)= \gamma +b'(u)\leq M, \s \s \mbox{for all $u \in R \,,$ }
\eeq 
where $M>0$ a constant.
\medskip

If $u(x) \in H^2(\Omega) \cap H^1_0(\Omega)$ is a solution of (\ref{2}), we decompose it as 
\beq
\lbl{5}
u(x)= \Sigma _{i=1}^n  \xi _i \p _i+U (x),
\eeq
where $U (x)$ is orthogonal to $\p _1, \ldots, \p _n$ in $L^2(\Omega)$. 
\medskip

For the problem (\ref{2}) we pose an inverse problem: keeping $e(x)$ fixed, find $\mu=\left( \mu _1, \ldots, \mu _n \right)$ so that the problem (\ref{2}) has a solution of any prescribed $n$-signature $\xi=\left( \xi _1, \ldots, \xi _n \right)$.

\begin{thm}\lbl{thm:1}
For the problem (\ref{2}) assume that the conditions (\ref{4}), (\ref{3}) hold, and 
\[
kM<\la _{n+1} \,.
\]
 Then given any $\xi=\left( \xi _1, \ldots, \xi _n \right)$, one can find a unique $\mu=\left( \mu _1, \ldots, \mu _n \right)$ for which the problem (\ref{2}) has a solution $u(x) \in H^2(\Omega) \cap H^1_0(\Omega)$ of $n$-signature $\xi$. This solution is unique. Moreover, we have a continuous curve of solutions $(u(k),\mu(k))$, such that  $u(k)$ has a fixed $n$-signature $\xi$, for all $0 \leq k \leq 1$.
\end{thm}

\pf
Let $e(x)=\Sigma _{j=n+1}^{\infty} e_j \p _j$. When $k=0$, the unique solution of (\ref{2}) of signature $\xi$ is $u(x)=\Sigma _{j=1}^{n} \xi _j \p _j-\Sigma _{j=n+1}^{\infty} \frac{e_j}{\la _j} \p _j$, corresponding to $\mu _j =-\la _j \xi _j$, $j=1, \ldots, n$. We shall use the implicit function theorem to continue this solution in $k$. With $u(x)= \Sigma _{i=1}^n  \xi _i \p _i+U (x)$, we multiply the equation (\ref{2}) by $\p _i$, and integrate
\beq
\lbl{16}
\mu _i=-\la _i \xi _i+k \int _\Omega g \left( \Sigma _{i=1}^n  \xi _i \p _i+U \right) \p _i \, dx, \s i=1, \ldots, n \,.
\eeq
Using these expressions in (\ref{2}), we have
\beq
\lbl{17}
\s \s \s \Delta U+kg\left( \sum _{i=1}^n  \xi _i \p _i+U \right)-k \Sigma _{i=1}^n \int _\Omega g \left( \sum _{i=1}^n  \xi _i \p _i+U \right) \p _i \, dx \p _i=e(x), 
\eeq
\[
 U=0 \s \mbox{on $\partial \Omega$} \,.
\]
The equations (\ref{16}) and (\ref{17}) constitute the classical Lyapunov-Schmidt decomposition of our problem (\ref{2}).
Define $H^2_{{\bf 0}}$ to be the subspace of $H^2(\Omega) \cap H^1_0(\Omega)$, consisting of functions with zero $n$-signature:
\[
H^2_{{\bf 0}} =  \left\{ u \in H^2(\Omega) \cap H^1_0(\Omega) \; | \; \int _\Omega u \p _i \, dx =0, \; i=1, \ldots, n \right\}.
\]
We recast the problem (\ref{17}) in the operator form as
\[
F(U, k) =e(x),
\]
where $ F(U, k) : H^2_{{\bf 0}}  \times R  \ra L^2(\Omega)$ is given by the left hand side of (\ref{17}). Compute the Frechet derivative
\[
F_{U}(U, k)w=\Delta w+kg' \left( \Sigma _{i=1}^n  \xi _i \p _i+U \right)w-\mu^*_1 \p _1-\ldots -\mu^* _n \p _n \,,
\]
where $\mu^*_i=k \int _\Omega g' \left( \Sigma _{i=1}^n  \xi _i \p _i+U \right) w \p _i \, dx$. By Lemma \ref{lma:3} the map $F_{U}(U, k)$ is injective. Since this map is Fredholm of index zero, it is also surjective. The implicit function theorem applies, giving us locally a curve of solutions $U=U(k)$. Then we compute $\mu=\mu (k)$ from (\ref{16}).
\medskip

To show that this curve can be continued for all $k$,  we only need to show that this curve $(u(k),\mu (k))$ cannot go to infinity at some $k$, i.e., we need an a priori estimate. Since the $n$-signature of the solution is fixed, we only need to estimate $U$. We claim that there is a constant $c>0$, so that 
\beq
\lbl{18}
||U||_{H^2(\Omega)} \leq c \,.
\eeq
We rewrite the equation in (\ref{17}) as
\[
\Delta U+k\gamma U= -kb\left( \sum _{i=1}^n  \xi _i \p _i+U \right)+k \Sigma _{i=1}^n \int _\Omega b \left( \sum _{i=1}^n  \xi _i \p _i+U \right) \p _i \, dx \p _i+e(x) \,.
\]
By the Corollary \ref{cor:2} to Lemma \ref{lma:3}, the estimate (\ref{18}) follows, since $b(u)$ is bounded.
\medskip

Finally, if the problem (\ref{2}) had a different solution $(\bar u(k),\bar \mu (k))$ with the same signature $\xi$, we would continue it back in $k$, obtaining at $k=0$ a different  solution of the linear problem of signature $\xi$ (since solution curves do not intersect by the implicit function theorem), which is impossible.
\epf

The Theorem \ref{thm:1} implies that the value of $\xi =(\xi_1, \ldots, \xi_n)$ uniquely identifies the solution pair $(\mu, u(x))$, where $\mu =(\mu_1, \ldots, \mu _n)$. Hence, the solution set of (\ref{2}) can be faithfully described by the map: $\xi \in R^n \ra \mu \in R^n$, which we call the {\em solution manifold}. In case $n=1$, we have the {\em solution curve} $\mu=\mu(\xi)$, which faithfully depicts the solution set. We show next that the solution manifold is connected.

\begin{thm}\lbl{thm:2}
In the conditions of Theorem \ref{thm:1}, the solution $(u,\mu_1,\dots,\mu _n)$ of (\ref{2}) is a continuous function of $\xi=(\xi_ 1, \dots ,\xi _n)$. Moreover,  we can continue solutions of any signature $\bar \xi$ to solution of arbitrary signature $\hat \xi $ by following any continuous curve in $R^n$ joining $\bar \xi$ and $\hat \xi$. 
\end{thm}

\pf
We use the implicit function theorem to show that any solution of (\ref{2}) can be continued in $\xi$. The proof is essentially the same as for continuation in $k$ above. After performing the same Lyapunov-Schmidt decomposition,
we recast the problem (\ref{17}) in the operator form
\[
F(U,\xi)=e(x) \,,
\]
where $F \, : \, H^2_{{\bf 0}} \times R^n  \ra L^2$ is defined  by the left hand side of (\ref{17}). The Frechet derivative $F_{U}(U, \xi)w$ is the same as before, and by the implicit function theorem we have locally $U=U(\xi)$. Then we compute $\mu=\mu (\xi)$ from (\ref{16}). We use the same a priori bound (\ref{18}) to continue the curve for all $\xi \in R^n$. (The bound (\ref{18}) is uniform in  $\xi$.)
\epf

Given a Fourier series $u(x)=\Sigma _{j=1}^{\infty} \xi _j \p _j$, we call the vector $(\xi _i, \ldots, \xi _n)$ to be the $(i,n)$-{\em signature } of $u(x)$.
Using Lemma \ref{lma:4} instead of Lemma \ref{lma:3}, we have the following variation of the above result.

\begin{thm}\lbl{thm:3}
For the problem (\ref{2}) assume that the conditions (\ref{4}), (\ref{3}) hold, and
\[
\la _{i-1} <k\gamma + kg'(u)< \la _{n+1}, \,  \s \mbox{for all $u \in R$} \,.
\]
 Then given any $\xi=\left( \xi _i, \ldots, \xi _n \right)$, one can find a unique $\mu=\left( \mu _i, \ldots, \mu _n \right)$ for which the problem
\beqa
\lbl{20}
& \Delta u+kg(u)=\mu_i \p _i+ \cdots +\mu _n \p _n+e(x), \, \s \mbox{for $x \in \Omega$}, \\ \nonumber
& u=0 \s \mbox{on $\partial \Omega$}
\eeqa
 has a solution $u(x) \in H^2(\Omega) \cap H^1_0(\Omega)$ of the $(i,n)$-signature $\xi$. This solution is unique. Moreover, we have a continuous curve of solutions $(u(k),\mu(k))$, such that  $u(k)$ has a fixed $(i,n)$-signature $\xi$, for all $0 \leq k \leq 1$. In addition, we can continue solutions of any  $(i,n)$-signature $\bar \xi$ to solution of arbitrary  $(i,n)$-signature $\hat \xi $ by following any continuous curve in $R^{n-i+1}$ joining $\bar \xi$ and $\hat \xi$.
\end{thm}

\section{Unbounded perturbations at resonance}
\setcounter{equation}{0}
\setcounter{thm}{0}
\setcounter{lma}{0}

We use an idea from \cite{I} to get the following a priori estimate.
\begin{lma}\lbl{lma:6}
Let $u(x)$ be a solution of the problem
\beq
\lbl{22}
\Delta u +\la _1 u+a(x)u=\mu _1 \p _1+e(x) \s \mbox{on $\Omega$}, \s u=0 \s \mbox{on $\partial \Omega$},
\eeq
with $e(x) \in \p _1 ^\perp$, and $a(x) \in C(\Omega)$. Assume  there is a constant $\gamma$, so that
\[
0 \leq a(x)\leq \gamma <\la_2 -\la _1, \s\s \mbox{for all $x \in \Omega$} \,.
\]
Write the solution of (\ref{22}) in the form $u(x)=\xi _1 \p _1+U$, with $U \in \p _1 ^\perp$, and assume that
\beq
\lbl{22a}
\xi _1 \mu _1 \leq 0 \,.
\eeq
Then there exists a constant $c_0$, so that
\beq
\lbl{22.2}
\int_\Omega |\nabla U|^2 \, dx \leq c_0 \,, \s\s \mbox{uniformly in $\xi _1 $ satisfying (\ref{22a})}\,.
\eeq
\end{lma}

\pf
We have
\beq
\lbl{22.1}
\s \s \Delta U +\la _1 U+a(x)\left(\xi _1 \p _1+U \right)=\mu _1 \p _1+e(x) \s \mbox{on $\Omega$}, \s U=0 \s \mbox{on $\partial \Omega$} \,.
\eeq
Multiply this by $\xi _1 \p _1-U $, and integrate
\[
\int_\Omega \left(|\nabla U|^2- \la _1 U^2 \right)\, dx+\int_\Omega a(x) \left(\xi _1^2 \p _1^2-U^2 \right) \, dx-\xi _1 \mu _1=-\int_\Omega eU \, dx \,.
\]
Dropping two non-negative terms on the left, we have
\[
\left(\la _2-\la _1-\gamma \right)\int_\Omega  U^2 \, dx \leq \int_\Omega \left(|\nabla U|^2- \la _1 U^2 \right)\, dx-\int_\Omega a(x) U^2 \, dx \leq -\int_\Omega eU \, dx \,.
\]
From this we get an estimate on $\int_\Omega  U^2 \, dx$, and then on $\int_\Omega |\nabla U|^2 \, dx$.
\epf

\begin{cor}\lbl{cor:4}
If, in addition, $\mu _1=0$ and $e(x) \equiv 0$, then $U \equiv 0$.
\end{cor}

We now consider the problem
\beq
\lbl{23}
\Delta u +\la _1 u+g(u)=\mu _1 \p _1+e(x) \s \mbox{on $\Omega$}, \s u=0 \s \mbox{on $\partial \Omega$} \,,
\eeq
with $e(x) \in \p _1 ^\perp$. We wish to find a solution pair $(u, \mu _1)$. We have the following extension of the result of R. Iannacci et al \cite{I}.

\begin{thm}\lbl{thm4}
Assume that $g(u) \in C^1(R)$ satisfies
\beq
\lbl{24}
u g(u) >0 \s\s \mbox{for all $u \in R$} \,,
\eeq
\beq
\lbl{25}
 g'(u) \leq \gamma< \la _2-\la _1 \s\s \mbox{for all $u \in R$} \,.
\eeq
Then there is a continuous curve of solutions of (\ref{23}): $(u(\xi _1),\mu _1(\xi _1))$, $u \in H^2(\Omega) \cap H^1_0(\Omega)$, with $-\infty<\xi _1<\infty$, and $\int _\Omega u(\xi _1) \p _1 \, dx=\xi _1$. This curve exhausts the solution set of (\ref{23}). The continuous function $\mu _1(\xi _1)$ is positive for $\xi _1 >0$ and large, and $ \mu _1(\xi _1)<0$ for $\xi _1 <0$ and $|\xi _1|$ large. In particular, $\mu _1(\xi^0 _1)=0$ at some $\xi^0 _1$, i.e., we have a solution of 
\[
\Delta u +\la _1 u+g(u)=e(x) \s \mbox{on $\Omega$}, \s u=0 \s \mbox{on $\partial \Omega$} \,.
\]
\end{thm}

\pf
By the Theorem \ref{thm:1} there exists a curve of solutions of (\ref{23}) $(u(\xi _1),\mu _1(\xi _1))$, which 
exhausts the solution set of (\ref{23}). The condition (\ref{24}) implies that $g(0)=0$, and then integrating (\ref{25}), we conclude that 
\beq
\lbl{27}
0 \leq \frac{g(u)}{u} \leq \gamma< \la _2-\la _1, \s\s \mbox{for all $u \in R$} \,. 
\eeq
Writing $u(x)=\xi _1 \p _1+U$, with $U \in \p _1 ^\perp$, we see that $U$ satisfies
\[
\Delta U +\la _1 U+g(\xi _1 \p _1+U)=\mu _1 \p _1+e(x) \s \mbox{on $\Omega$}, \s U=0 \s \mbox{on $\partial \Omega$} \,.
\]
We rewrite this equation in the form (\ref{22}), by letting $a(x)=\frac{ g(\xi _1 \p _1+U)}{\xi _1 \p _1+U}$. By (\ref{27}), the Lemma \ref{lma:6} applies, giving us the estimate (\ref{22.2}).
\medskip

We claim next that $|\mu _1(\xi _1)|$ is bounded uniformly in $\xi _1$, provided that $\xi _1 \mu _1 \leq 0$. Indeed, let us assume first that $\xi _1 \geq 0$ and $\mu _1 \leq 0$. Then
\[
\mu_1=\int_\Omega g(u) \p _1 \, dx = \int_\Omega \frac{g(u)}{u} \xi _1  \p _1^2 \, dx+\int_\Omega \frac{g(u)}{u} U \p _1  \, dx \geq  \int_\Omega \frac{g(u)}{u} U \p _1  \, dx\,,
\]
\[
|\mu _1|=-\mu _1 \leq -\int_\Omega \frac{g(u)}{u} U \p _1  \, dx \leq \gamma \int_\Omega | U \p _1 | \, dx\leq c_1 \,,
\]
for some $c_1>0$, in view of (\ref{27}) and  the estimate (\ref{22.2}). The case when $\xi _1 \leq 0$ and $\mu _1 \geq 0$ is similar.
\medskip

We now rewrite (\ref{23}) in the form
\beq
\lbl{28}
\Delta u+a(x)u=f(x)  \s \mbox{on $\Omega$}, \s u=0 \s \mbox{on $\partial \Omega$}\,
\eeq
with $a(x)=\la _1 +\frac{g(u)}{u}$, and $f(x)=\mu _1 \p _1+e(x)$. By above, we have a uniform in $\xi _1$ bound on $||f||_{L^2(\Omega)}$, and by the Corollary \ref{cor:4} we have uniqueness for (\ref{28}). It follows that
\[
||u||_{H^2(\Omega)} \leq c||f||_{L^2(\Omega)} \leq c_2 \,,
\]
for some $c_2>0$.
\medskip

Assume, contrary to what we wish to prove, that there is a sequence $\{\xi _1^n \} \ra \infty$, such that $\mu _1 (\xi _1^n) \leq 0$. We have 
\[
u=\xi^n _1 \p _1+U \,,
\]
with both $u$ and $U$ bounded in $L^2(\Omega)$, uniformly in $\xi _1^n$,
which results in a contradiction for $n$ large. We prove similarly that $ \mu _1(\xi _1)<0$ for $\xi _1 <0$ and $|\xi _1|$ large.
\epf

\noindent
{\bf Example} We have solved numerically  the problem 
\beqa \nonumber
& u''+u+0.2 \, \frac{u^3}{u^2+3u+3}+\sin \frac12 u=\mu \sin x+ 5 \left(x-\pi/2 \right), \s 0<x<\pi, \\ \nonumber
& u(0)=u(\pi)=0 \,. \nonumber
\eeqa
The Theorem \ref{thm4} applies. Write the solution as $u(x)=\xi \sin x +U(x)$, with $\int_0^{\pi} U(x) \sin x \, dx=0$. 
Then the solution curve  $\mu=\mu (\xi)$ is given in Figure $1$. The picture suggests that the problem has at least one  solution for all $\mu$.

\begin{figure}
\begin{center}
\scalebox{0.9}{\includegraphics{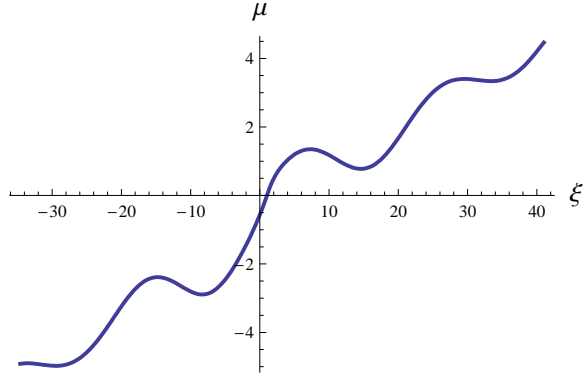}}
\end{center}
\caption{ An example for the Theorem \ref{thm4}}
\end{figure}

\medskip

We have the following extension of the results of D.G. de Figueiredo and W.-M.  Ni \cite{FN} and R. Iannacci et al \cite{I}, which does not require that $\mu =0$.

\begin{thm}\lbl{thm:8}
In addition to the conditions of the Theorem \ref{thm4}, assume that for some constants  $c_0>0$ and $p> \frac32$
\beq
\lbl{28.1}
ug(u) > c_0 |u|^p, \s \mbox{for all $u >0$ ($u <0$)} \,.
\eeq
Then for the problem  (\ref{23}) we have  $\lim _{\xi _1 \ra \infty} \mu (\xi _1)=\infty$ ($\, \lim _{\xi _1 \ra -\infty} \mu (\xi _1)=-\infty$).
\end{thm}

\pf
Assume that (\ref{28.1}) holds for $u >0$. By the Theorem \ref{thm4}, $ \mu (\xi _1)>0$ for $\xi _1$ large. Assume, on the contrary, that $\mu (\xi _1)$ is bounded along some sequence of $\xi _1$'s, which tends to $\infty$. Writing $u=\xi _1 \p _1+U$, we conclude from the line following (\ref{22.1}) that
\beq
\lbl{28.2}
\int_\Omega U^2  \, dx \leq c_1 \xi _1+c_2, \s \mbox{for some constants  $c_1>0$ and $c_2>0$} \,.
\eeq
We have
\[
\mu_1=\int_\Omega g(\xi _1 \p _1+U) \p _1 \, dx =\int_\Omega \left( g(\xi _1 \p _1+U)-g(\xi _1 \p _1) \right) \p _1 \, dx+\int_\Omega g(\xi _1 \p _1) \p _1 \, dx  \,.
\]
Using the mean value theorem, the estimate (\ref{28.2}), and the condition (\ref{28.1}), we estimate
\[
\mu_1 >c_3 {\xi _1}^{p-1}-c_4 {\xi _1}^{1/2}-c_5 \,,
\]
with some positive constants $c_3$, $c_4$ and $c_5$. It follows that $ \mu (\xi _1)$ gets large along our sequence, a contradiction.
\epf

Bounded perturbations at resonance are much easier to handle. For example, we have the following result.

\begin{thm}\lbl{thm:5}
Assume that $g(u) \in C^1(R)$ is a bounded function, which  satisfies the condition (\ref{24}), and in addition,
\[
 \lim _{u \ra \pm \infty} g(u)=0 \,.
\]
There is a continuous curve of solutions of (\ref{23}): $(u(\xi _1),\mu _1(\xi _1))$, $u \in H^2(\Omega) \cap H^1_0(\Omega)$, with $-\infty<\xi _1<\infty$, and $\int _\Omega u(\xi _1) \p _1 \, dx=\xi _1$. This curve exhausts the solution set of (\ref{23}).
Moreover,  there are constants $\mu _- <0< \mu _+$ so that the problem (\ref{23}) has at least two   solutions for $\mu \in (\mu _-,\mu _+) \setminus 0$, it has at least one   solution for $\mu=\mu _- $, $\mu=0$ and $\mu=\mu _+ $, and no  solutions for $\mu$ lying outside of $(\mu _-,\mu _+)$.
\end{thm}

\pf
Follow the proof of the Theorem \ref{thm4}. Since $g(u)$ is bounded, we have a uniform in $\xi _1$ bound on $||U||_{C^1}$, see \cite{FN}.  Since $\mu_1=\int_\Omega g(\xi _1 \p _1+U) \p _1 \, dx$, we conclude that for $\xi _1$ positive (negative) and large, $\mu _1$ is positive (negative) and it tends to zero as $\xi _1 \ra \infty$ ($\xi _1 \ra -\infty$).
\epf

\noindent
{\bf Example} We have solved numerically  the problem 
\[
 u''+u+\frac{u}{2u^2+u+1}=\mu \sin x+ \sin 2x, \s 0<x<\pi, 
\s u(0)=u(\pi)=0 \,.
\]
The Theorem \ref{thm:5} applies. Write the solution as $u(x)=\xi \sin x +U(x)$, with $\int_0^{\pi} U(x) \sin x \, dx=0$. 
Then the solution curve  $\mu=\mu (\xi)$ is given in Figure $2$. The picture shows that, say, for $\mu=-0.4$, the problem has exactly two solutions, while for $\mu=1$ there are no solutions.

\begin{figure}
\begin{center}
\scalebox{0.9}{\includegraphics{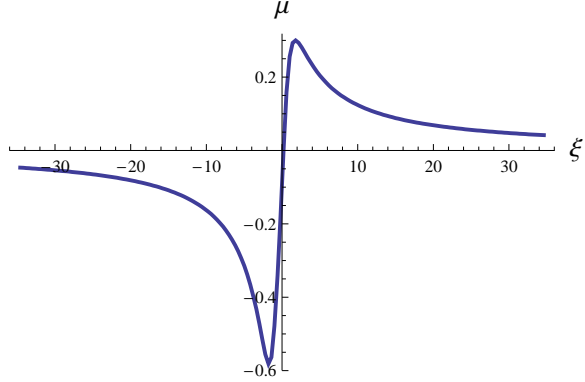}}
\end{center}
\caption{ An example for the Theorem \ref{thm:5}}
\end{figure}

\medskip

We also have  a result of Landesman-Lazer type, which also provides some additional information on the solution curve.

\begin{thm}\lbl{thm:12}
Assume that the function $g(u) \in C^1(R)$ is bounded, it satisfies  (\ref{25}), and in addition,  $g(u)$ has finite limits at $\pm \infty$, and
\[
g(-\infty)<g(u)<g(\infty), \, \s \mbox{for all $u \in R$} \,.
\]
Then there is a continuous curve of solutions of (\ref{23}): $(u(\xi _1),\mu _1(\xi _1))$, $u \in H^2(\Omega) \cap H^1_0(\Omega)$, with $-\infty<\xi _1<\infty$, and $\int _\Omega u(\xi _1) \p _1 \, dx=\xi _1$. This curve exhausts the solution set of (\ref{23}), and $\lim _{\xi _1 \ra \pm \infty}\mu _1(\xi _1) = g(\pm \infty) \int _\Omega \p _1 \, dx$. I.e.,  the problem (\ref{23}) has a  solution if and only if 
\[
g(- \infty) \int _\Omega \p _1 \, dx<\mu<g( \infty) \int _\Omega \p _1 \, dx \,.
\]
\end{thm}

\pf
Follow the proof of the Theorem \ref{thm4}. Since $g(u)$ is bounded,  we have a uniform bound on $U$, when we do the continuation in $\xi _1$. Hence $\mu _1 \ra g(\pm \infty) \int _\Omega \p _1 \, dx$, as $\xi _1 \ra \pm \infty$, and by continuity of $\mu _1 (\xi _1)$, the problem (\ref{23}) is solvable for all $\mu _1$'s lying between these limits.
\epf

\noindent
{\bf Example} We have solved numerically  the problem 
\[
 u''+u+\frac{u}{\sqrt{u^2+1}}=\mu \sin x+ 5\sin 2x-\sin 10 x, \s 0<x<\pi, 
\s u(0)=u(\pi)=0 \,.
\]
The Theorem \ref{thm:12} applies. Write the solution as $u(x)=\xi \sin x +U(x)$, with $\int_0^{\pi} U(x) \sin x \, dx=0$. Then the solution curve  $\mu=\mu (\xi)$ is given in Figure $3$. It confirms that $\lim _{\xi _1 \ra \pm \infty}\mu _1(\xi _1)  =\pm \frac{4}{\pi}$ ($\frac{4}{\pi}=\int _0^{\pi} \frac{2}{\pi} \, \sin x \, dx$).

\begin{figure}
\begin{center}
\scalebox{0.9}{\includegraphics{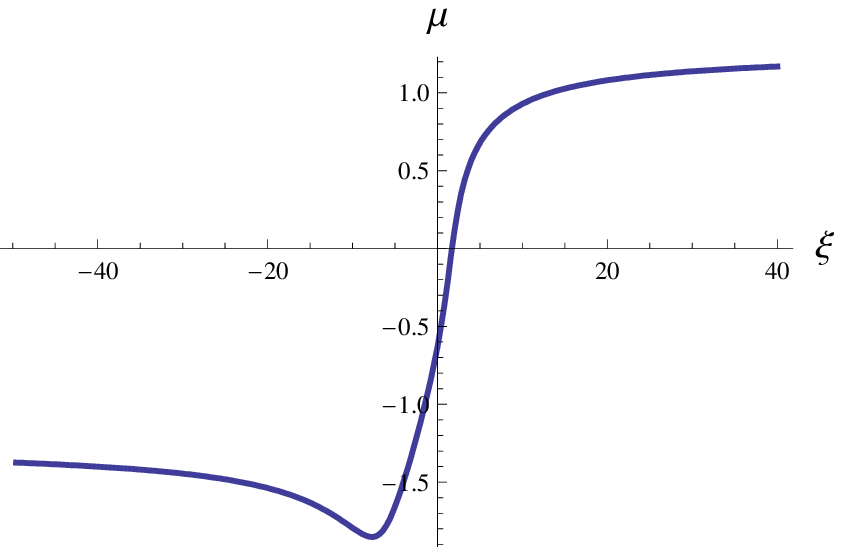}}
\end{center}
\caption{ An example for the Theorem \ref{thm:12}}
\end{figure}

One can append the following uniqueness condition (\ref{28.50}) to all of the above results. For example, we have the following result.

\begin{thm}
Assume that the conditions of the Theorem \ref{thm4} hold, and in addition 
\beq
\lbl{28.50}
g'(u)>0, \, \s \mbox{for all $u \in R$} \,.
\eeq
Then
\beq
\lbl{28.5}
\mu _1'(\xi _1)>0, \s \mbox{ for all $\xi _1 \in R$} \,.
\eeq
\end{thm}

\pf
Clearly, $\mu' _1(\xi _1)>0$ at least for some values of $\xi _1$. If (\ref{28.5}) is not true, then  $\mu' _1(\xi^0 _1)=0$ at some $\xi^0 _1$.
Differentiate the equation (\ref{23}) in $\xi _1$, set $\xi _1=\xi^0 _1$, and denote $w=u_{\xi _1} |_{\xi _1=\xi^0 _1}$, obtaining
\beqa \nonumber
& \Delta w+\left(\la _1  +g'(u) \right)w=0  \s \mbox{for $x \in \Omega$}, \\ \nonumber
& w=0 \s \mbox{on $\partial \Omega$}.
\eeqa
Clearly, $w$ is not zero, since it has a non-zero projection on $\p _1$ ($U_{\xi _1} \in \p _1^{\perp}$). On the other hand, $w \equiv 0$, since by the assumption (\ref{25}) we have $\la _1<\la _1  +g'(u)<\la _2$.
\epf

\begin{cor}
In addition to the conditions of this theorem, assume that the condition (\ref{28.1}) holds, for all $u \in R$. Then for any $f(x) \in L^2(\Omega)$, the problem
\[
\Delta u+\la _1 u +g(u) =f(x) \,,  \s \mbox{for $x \in \Omega$}, \s u=0 \s \mbox{on $\partial \Omega$}
\]
has a unique solution $u(x) \in H^2(\Omega) \cap H^1_0(\Omega)$.
\end{cor}

\section{Resonance at higher eigenvalues}
\setcounter{equation}{0}
\setcounter{thm}{0}
\setcounter{lma}{0}

We  consider the problem
\beq
\lbl{43}
\Delta u +\la _k u+g(u)=f(x) \s \mbox{on $\Omega$}, \s u=0 \s \mbox{on $\partial \Omega$} \,,
\eeq
where $\la _k$ is assumed to be a {\em simple} eigenvalue of $-\Delta$. We have the following extension of the result of D.G. de Figueiredo and W.-M.  Ni \cite{FN} to the case of resonance at a non-principal eigenvalue.

\begin{thm}\lbl{thm7}
Assume that $g(u) \in C^1(R)$ is  bounded, it satisfies (\ref{24}), and
\beq
\lbl{44}
g'(u) \leq c_0, \s \mbox{for all $u \in R$, and some $c_0>0$} \,,
\eeq
\beq
\lbl{45}
\liminf _{u \ra \infty} g(u)>0, \s \limsup _{u \ra -\infty} g(u)<0 \,.
\eeq
Assume that $f(x) \in L^2(\Omega)$ satisfies
\beq
\lbl{46}
\int _\Omega f(x) \p _k (x) \, dx=0 \,.
\eeq
Then the problem (\ref{43}) has a solution $u(x) \in H^2(\Omega) \cap H^1_0(\Omega)$.
\end{thm}

\pf
By (\ref{44}) we may assume that $\la _k+ g'(u) \leq \la _{n+1}$ for some $n>k$. Expand $f(x)=\mu _1^0 \p _1+\mu _2^0 \p _2+ \cdots +\mu _n^0 \p _n+e(x)$, with 
$e(x) \in Span \{ \p_1, \ldots, \p _{n} \}^{\perp}$, and $u(x)=\xi _1 \p _1+\xi _2 \p _2+ \cdots +\xi _n \p _n+U(x)$, with 
$U(x) \in Span \{ \p_1, \ldots, \p _{n} \}^{\perp}$. By (\ref{46}), $\mu ^0 _k=0$. By the Theorem \ref{thm:1} for any $\xi=\left( \xi _1, \ldots, \xi _n \right)$, one can find a unique $\mu=\left( \mu _1, \ldots, \mu _n \right)$ for which the problem (\ref{2}) has a solution of $n$-signature $\xi$, and we need to find a $\xi ^0=\left( \xi^0 _1, \ldots, \xi^0 _n \right)$, for which $\mu (\xi ^0)=\left( \mu^0 _1, \ldots, \mu^0 _{k-1},0,\mu^0 _{k+1}, \ldots, \mu^0 _n \right)$.
\medskip

Multiplying the equation  (\ref{43}) by $\p _i$, and integrating we get 
\[
(\la _k-\la _i) \xi _i+ \int _\Omega g \left(\sum _{i=1}^n \xi _i \p _i +U\right) \p _i \, dx=\mu ^0 _i, \s i=1, \ldots, k-1,k+1, \ldots n 
\]
\[
\int _\Omega g \left(\sum _{i=1}^n \xi _i \p _i +U\right) \p _k\, dx=0 \,.
\]
We need to solve this system of equations for $\left( \xi _1, \ldots, \xi _n \right)$. For that we set up a map $T : \left( \eta _1, \ldots, \eta _n \right) \ra \left( \xi _1, \ldots, \xi _n \right)$, by calculating $\xi _i$ from
\[
(\la _k-\la _i) \xi _i= \mu ^0 _i -\int _\Omega g \left(\sum _{i=1}^n \eta _i \p _i +U\right) \p _i\, dx, \s i=1, \ldots, k-1,k+1, \ldots n
\]
followed by
\[
\xi _k=\eta _k-\int _\Omega g \left(\xi _1 \p _1+ \cdots +\xi _{k-1} \p _{k-1}+\eta _k \p _k +\xi _{k+1} \p _{k+1}+ \cdots +\xi _n \p _n+U \right)\p _k \, dx \,.
\]
Fixed points of this map provide  solutions to our system of equations. By the Theorem \ref{thm:2}, the map $T$ is continuous. Since $g(u)$ is bounded, $\left( \xi _1, \ldots,\xi _{k-1},\xi_{k+1},\ldots, \xi _n \right)$ belongs to a bounded set. By (\ref{24}) and (\ref{45}), $\xi _k <\eta _k$ for $\eta _k>0$ and large, while $\xi _k >\eta _k$ for $\eta _k<0$ and $|\eta _k|$ large. Hence, the map $T$ maps a sufficiently large ball around the origin in $R^n$ into itself, and Brouwer's fixed point theorem applies, giving us a fixed point of $T$.
\epf

\section{Numerical computation of solutions}
\setcounter{equation}{0}
\setcounter{thm}{0}
\setcounter{lma}{0}

We describe numerical computation of solutions for the problem
\beq
\lbl{n1}
u''+u+g(u)=\mu \sin x+e(x), \s 0<x<\pi, \s u(0)=u(\pi)=0 \,,
\eeq
whose linear part is at resonance. We assume that $\int _0^{\pi} e(x) \sin x \, dx=0$. Writing $u(x)=\xi \sin x +U(x)$, with $\int _0^{\pi} U(x) \sin x \, dx=0$, we shall compute the solution curve of (\ref{n1}): $(u(\xi),\mu (\xi))$. (I.e., we write $\xi$, $\mu$ instead of $\xi _1$, $\mu _1$.) We shall use Newton's method to perform continuation in $\xi$.
\medskip

Our first task is to implement the ``linear solver", i.e., the numerical solution of the following problem: given any $\xi \in R$, and any functions $a(x)$ and $f(x)$, find $u(x)$ and $\mu$ solving
\beqa
\lbl{n2}
& u''+a(x)u=\mu \sin x+f(x), \s 0<x<\pi \,, \\\nonumber
& u(0)=u(\pi)=0 \,,\\ \nonumber
& \int _0^{\pi} u(x) \sin x \, dx=\xi \,.\nonumber
\eeqa
The general solution of the equation (\ref{n2}) is of course
\[
u(x)=Y(x)+c_1 u_1(x)+c_2 u_2(x) \,,
\]
where $Y(x)$ is any particular solution, and $u_1$,$u_2$ are two solutions of the corresponding homogeneous equation
\beq
\lbl{n3}
u''+a(x)u=0, \s 0<x<\pi \,.
\eeq
We shall use $Y=\mu Y_1+ Y_2$, where $Y_1$ solves 
\[
u''+a(x)u=\sin x, \s u(0)=0, \s u'(0)=1 \,,
\]
and   $Y_2$ solves
\[
u''+a(x)u=f(x), \s u(0)=0, \s u'(0)=1 \,.
\]
Let $u_1(x)$ be the solution of (\ref{n3}) with  $u(0)=0$, $u'(0)=1$, and let  $u_2(x)$ be any solution of (\ref{n3}) with  $u_2(0) \ne 0$.
The condition $u(0)=0$ implies that $c_2=0$, i.e., there is no need to compute $u_2(x)$, and we have  
\beq
\lbl{n4}
u(x)=\mu Y_1(x)+Y_2(x)+c_1 u_1(x) \,.
\eeq
We used the NDSolve command in {\em Mathematica} to calculate $u_1$, $Y_1$ and $Y_2$. {\em Mathematica} not only solves differential equations numerically, but it returns the solution as an interpolated function of $x$, practically indistinguishable from an explicitly defined function.The condition $u(\pi)=0$ and the last line in (\ref{n2}) imply that
\[
\mu Y_1(\pi)+c_1 u_1(\pi)=-Y_2(\pi) \,,
\]
\[
\mu \int _0^{\pi} Y_1(x)\sin x \, dx+c_1 \int _0^{\pi} u_1(x)\sin x \, dx=\xi -\int _0^{\pi} Y_2(x)\sin x \, dx \,,
\]
Solving this system for $\mu$ and $c_1$, and using them in (\ref{n4}), we obtain the solution of (\ref{n2}).
\medskip

Turning to the problem (\ref{n1}), we begin with an initial $\xi _0$, and using a step size $\Delta \xi$, on a mesh $\xi _i=\xi _0 +i \Delta \xi$, $i=1,2, \ldots, nsteps$,  we compute the solution of (\ref{n1}), satisfying $\int _0^{\pi} u (x) \sin x \, dx=\xi _i$, by using Newton's method. Namely, assuming that the iterate $u_n(x)$ is already computed, we linearize the equation (\ref{n1}) at it, i.e., we solve the problem (\ref{n2}) with
$a(x)=1+g'(u_n(x))$, $f(x)=-g(u_n(x))+g'(u_n(x)) u_n(x)+e(x)$, and $\xi =\xi _i$. After several iterations, we compute $(u(\xi _i), \mu (\xi _i))$. We found that two iterations of Newton's method, coupled with   $\Delta \xi$  not too large (e.g., $\Delta \xi=0.5$), were sufficient for accurate computation of the solution curves. To start Newton's iterations, we used $u(x)$ computed at the preceding step, i.e., $u(\xi _{i-1})$.
\medskip

We have verified our numerical results by an independent calculation. Once  a solution of (\ref{n1}) was computed at some $\xi _i$, we took its initial data $u(0)=0$ and $u'(0)$, and computed numerically the solution of the equation in (\ref{n1}) with this initial data, let us call it $v(x)$ (using the NDSolve command). We always had  $v(\pi)=0$ and $\int _0^{\pi} v(x) \sin x \, dx=\xi _i$.

\end{document}